%% file: dual.tex
\documentclass[11pt]{article}
\input{format} \input{mathdefs} \input{thmdefs}

\title{Amenability for dual Banach algebras}
\author{\it Volker Runde}
\date{}
\begin{document}
\maketitle
\begin{abstract}
We define a Banach algebra ${\mathfrak A}$ to be dual if ${\mathfrak A}
= ({\mathfrak A}_\ast )^\ast$ for a closed submodule ${\mathfrak A}_\ast$ of 
${\mathfrak A}^\ast$. The class of dual Banach algebras includes all 
$W^\ast$-algebras, but also all algebras $M(G)$ for locally compact groups 
$G$, all algebras ${\cal L}(E)$ for reflexive Banach spaces $E$, as well as 
all biduals of Arens regular Banach algebras. The general impression is that 
amenable, dual Banach algebras are rather the exception than the rule. 
We confirm this impression. We first show that under certain conditions an amenable
dual Banach algebra is already super-amenable and thus finite-dimensional. We then develop two notions of
amenability --- Connes-amenability and strong Connes-amenability --- which take the 
$w^\ast$-topology on dual Banach algebras into account. We relate the amenability of
an Arens regular Banach algebra $\A$ to the (strong) Connes-amenability of $\A^{\ast\ast}$;
as an application, we show that there are reflexive Banach spaces with the approximation
property such that ${\cal L}(E)$ is not Connes-amenable. We characterize the amenability of
inner amenable locally compact groups in terms of their algebras of pseudo-measures. Finally, 
we give a proof of the known fact that the amenable von Neumann algebras are the subhomogeneous ones which 
avoids the equivalence of amenability and nuclearity for $\cstar$-algebras. 
\end{abstract}
\begin{classification}
{\!}43A10, 46B28, 46H25, 46H99 (primary), 46L10, 46M20.
\end{classification}
\section{Introduction}
Amenable Banach algebras were introduced by B.\ E.\ Johnson in \cite{Joh1}, and
have since then turned out to be extremely interesting objects of research. 
The definition of an amenable Banach algebra is strong enough to allow for the development of a
rich general theory, but still weak enough to include a variety of interesting
examples. Very often, for a class of Banach algebras, the amenability condition
singles out an important subclass: For a locally compact group $G$, the 
convolution algebra $L^1(G)$ is amenable if and only if $G$ is amenable in
the classical sense (\cite{Joh1}); a $\cstar$-algebra is amenable if and
only if it is nuclear (\cite{Con}, \cite{BP}, \cite{Haa}); and a uniform algebra with
character space $\Omega$ is amenable if and only if it is 
${\cal C}_0(\Omega)$ (\cite{She}). To determine, for a given
class of Banach algebras, which algebras in it are the amenable ones is an
active area of research: For instance, it is still open for which Banach
spaces $E$ the Banach algebra ${\cal K}(E)$ of all compact operators on $E$
is amenable (see \cite{GJW} for partial results), and only recently, 
examples of radical amenable Banach algebras were given (\cite{Run}, \cite{Rea}).
\par
In this paper, we consider the following class of Banach algebras:
\begin{definition} \label{ddef}
A Banach algebra $\mathfrak A$ is said to be {\it dual\/} if there is a closed 
submodule ${\mathfrak A}_\ast$ of ${\mathfrak A}^\ast$ such that
${\mathfrak A} = ({\mathfrak A}_\ast)^\ast$.
\end{definition}
\par
If ${\mathfrak A}$ is a dual Banach algebra, the predual module 
${\mathfrak A}_\ast$ need not be unique. In this paper, however, it is
always clear, for a dual Banach algebra ${\mathfrak A}$, to which 
${\mathfrak A}_\ast$ we are referring. In particular, we may speak of
the $w^\ast$-topology on ${\mathfrak A}$ without ambiguity.
\par
The notion of a dual Banach algebra as defined in Definition \ref{ddef} 
is by no means universally accepted. The name ``dual Banach algebra'' occurs 
in the literature in several contexts --- often quite far apart from
Definition \ref{ddef}. On the other hand, Banach algebras satisfying
Definition \ref{ddef} may appear with a different name tag; for instance,
dual Banach algebras in our sense are called Banach algebras with $(DM)$
in \cite{CG2} and \cite{CG1}. 
\par
We note a few elementary properties of dual Banach algebras:
\begin{proposition}
Let $\mathfrak A$ be a dual Banach algebra. Then:
\begin{items}
\item Multiplication in $\mathfrak A$ is separately $w^\ast$-continuous.
\item $\mathfrak A$ has an identity if and only if it has a bounded approximate
identity.
\item The Dixmier projection $\pi \!: {\mathfrak A}^{\ast\ast} \cong
{\mathfrak A}_\ast^{\ast\ast\ast} \to {\mathfrak A}^\ast_\ast \cong
{\mathfrak A}$ is an algebra homomorphism with respect to either Arens
multiplication on ${\mathfrak A}^{\ast\ast}$.
\end{items}
\end{proposition}
\begin{proof}
(i) and (ii) are obvious, and (iii) follows from (i) and 
\cite[Theorem 1]{Pal}.
\end{proof}
\par
The main reason for us to consider dual Banach algebras is that this class
covers a wide range of examples:
\begin{examples}
\item Any $W^\ast$-algebra is dual.
\item If $G$ is a locally compact group, then $M(G)$ is dual (with
$M(G)_\ast = {\cal C}_0(G)$).
\item If $E$ is a reflexive Banach space, then ${\cal L}(E)$ is dual
(with ${\cal L}(E)_\ast = E \Tensor E^\ast$).
\item If $\mathfrak A$ is an Arens regular Banach algebra, then 
${\mathfrak A}^{\ast\ast}$ is dual; in particular, every reflexive
Banach algebra is dual.
\end{examples}
\par
Comparing this list of dual Banach algebras with our stock of amenable
Banach algebras, the overlap is surprisingly small. But although there are
few examples of dual Banach algebras which are known to be amenable,
there are equally few dual Banach algebras for which we positively know that
they are not amenable:
\begin{itemize}
\item A $W^\ast$-algebra is amenable if and only if it is subhomogeneous
(this follows from \cite[(1.9) Corollary]{Was} and \cite{Con}; see also
\cite{LLW}). Even for such a ``simple'' object as ${\mathbb M}_\infty :=
\mbox{$\ell^\infty$-}\bigoplus_{n=1}^\infty {\mathbb M}_n$, the proof of
non-amenability requires that amenability implies nuclearity.
\item If $G$ is a locally compact group, then $M(G)$ is amenable if and only 
if $G$ is discrete and amenable (\cite{DGH}).
\item The only Banach spaces $E$ for which ${\cal L}(E)$ is known to be 
amenable are the finite-dimensional ones, and they may well be the only ones. 
For a Hilbert space $\mathfrak H$, the results on amenable von Neumann algebras 
imply that ${\cal L}({\mathfrak H})$ is not amenable unless $\mathfrak H$ is finite-dimensional. 
It seems to be unknown, however, if ${\cal L}(\ell^p)$ is non-amenable for $p \in
(1,\infty) \setminus \{ 2 \}$.
\item The only known Arens regular Banach algebras $\mathfrak A$ for which
${\mathfrak A}^{\ast\ast}$ is amenable are the subhomogeneous 
$\cstar$-algebras; in particular, no infinite-dimensional, reflexive, amenable
Banach algebra is known. (It seems that the demand that ${\mathfrak A}^{\ast\ast}$
be amenable is very strong: It necessitates $\mathfrak A$ to be amenable
(\cite{GLW}, \cite{Gou}) and, for many classes of Banach algebras,
forces $\mathfrak A$ to be finite dimensional (\cite{GRW}, \cite{GLW},
\cite{Run}).)  
\end{itemize}
\par
The general impression thus is that amenability in the sense of \cite{Joh1} is 
too strong to
allow for the development of a rich theory for dual Banach algebras, and that
some notion of amenability taking the $w^\ast$-topology on dual Banach algebras
into account is more appropriate (\cite[Question 10]{Grb2}). 
Nevertheless, although amenability {\it seems\/} to be a condition which
is in conflict with Definition \ref{ddef}, this impression is supported by
surprisingly few proofs, and even where such proofs exist --- in the 
$W^\ast$-case, for instance ---, they often seem inappropriately deep. 
\par
This paper therfore aims into two directions: 
First, we want to substantiate our impression that dual Banach algebras are
rarely amenable with theorems, and secondly, we want to develop a suitable
notion of amenability --- which we shall call Connes-amenability --- for dual
Banach algebras.
\section{Amenability preliminaries} 
Let $\mathfrak A$ be a Banach algebra, and let $E$ be a Banach 
$\mathfrak A$-bimodule. A {\it derivation\/} from $\A$ into $E$ is a bounded linear map
satisfying
\[
  D(ab) = a \cdot Db + (Da) \cdot b \qquad (a,b \in \A).
\]
As is customary, we write ${\cal Z}^1(\A,E)$ for the Banach space of all derivations from
$\A$ into $E$. For $x \in E$, the linear map
\[
  \ad_x \!: \A \to E, \quad a \mapsto a \cdot x - x \cdot a
\]
is a derivation. Derivations of this form are called {\it inner derivations\/};
the normed space of all inner derivations from $\A$ to $E$ is denoted by
${\cal B}^1(\A,E)$. The quotient space ${\cal H}^1(\A,E) := {\cal Z}^1(\A,E)/
{\cal B}^1(\A,E)$ is called the {\it first cohomology group\/} (or rather: {\it space\/}) 
of $\A$ with coeffiecients in $E$.
\par
The dual space of a Banach $\A$-bimodule can be made into a Banach $\A$-bimodule as well via
\[
  \langle x, a \cdot \phi \rangle := \langle x \cdot a , \phi \rangle
  \quad\text{and}\quad
  \langle x, \phi \cdot a \rangle := \langle a \cdot x, \phi \rangle
  \qquad (a \in \A, \, \phi \in E^\ast, \, x \in E).
\]
A Banach algebra $\A$ is defined to be {\it amenable\/} if 
${\cal H}^1(\A, E^\ast) = \{ 0 \}$ for every Banach $\A$-bimodule $E$
(\cite{Joh1}). References for amenable Banach algebras
are \cite{Joh1} as well as \cite{Hel}, where a different, but equivalent approach,
based on the notion of flatness in topological homology is given.
\par
We shall also require a characterization of amenable Banach algebras in
terms of approximate diagonals as given in \cite{Joh2}.
Let ${\mathfrak A} \Tensor {\mathfrak A}$ denote the projective tensor product
of $\mathfrak A$ with itself. Then ${\mathfrak A} \Tensor {\mathfrak A}$
is a Banach $\mathfrak A$-bimodule through
\[
  a\cdot (x \tensor y) := ax \tensor y \quad\text{and}\quad \text
  (x \tensor y)\cdot a := x \tensor ya \qquad (a,x,y \in {\mathfrak A}).
\]
Let $\Delta \!: {\mathfrak A} \Tensor {\mathfrak A} \to {\mathfrak A}$ be
the multiplication operator, i.e.\ $\Delta(a \tensor b) := ab$ for $a,b 
\in {\mathfrak A}$ (sometimes, when we wish to emphasize the algebra
$\mathfrak A$, we also write $\Delta_{\mathfrak A}$). An 
{\it approximate diagonal\/} for $\mathfrak A$ is a bounded net 
$( m_\alpha )_\alpha$ in ${\mathfrak A} \Tensor {\mathfrak A}$ such that
\[
  a \cdot m_\alpha - m_\alpha \cdot a \to 0 \quad\text{and}\quad
  a\Delta m_\alpha \to a \qquad (a \in {\mathfrak A}).
\]
The algebra $\mathfrak A$ is amenable if and only if it has an approximate
diagonal (\cite{Joh2}). 
\par
There are several variants of amenability, two  of which we will discuss 
here: super-amenbility, and Connes-amenability.
\par
A Banach algebra $\mathfrak A$ is said to be {\it super-amenable\/} (or 
{\it contractible\/}) if ${\cal H}^1({\mathfrak A},E) = \{ 0 \}$ for
every Banach $\mathfrak A$-bimodule $E$. Equivalently, $\mathfrak A$ is 
super-amenable if it has a {\it diagonal\/}, i.e.\ a constant approximate
diagonal (\cite[Theorem 6.1]{CL}). All algebras ${\mathbb M}_n$ with
$n \in \posints$ and all finite direct sums of such algebras are 
super-amenable;
no other examples are known. Conversely, it is known that every super-amenable 
Banach algebra $\mathfrak A$ which satisfies some rather mild hypotheses in
terms of Banach space geometry must be a finite direct sum of full matrix
algebras (\cite{Sel}, for example; see \cite{Run} for a survey and some 
refinements). In particular, every super-amenable Banach algebra $\mathfrak A$
with the approximation property is of the form
\[
  {\mathfrak A} \cong {\mathbb M}_{n_1} \oplus \cdots \oplus
  {\mathbb M}_{n_k}
\]
with $n_1, \ldots, n_k \in \posints$.
\par
In \cite{JKR}, B.\ E.\ Johnson, R.\ V.\ Kadison, and J.\ Ringrose introduced a notion 
of amenability for von Neumann algebras which takes the ultraweak topology into account. 
The basic concepts, however, make sense for arbitrary dual Banach algebras. Let $\A$ be a dual Banach algebra, and
let $E$ be a Banach $\A$-bimodule. Then we call $E^\ast$ a {\it $w^\ast$-Banach
$\A$-bimodule\/} if, for each $\phi \in E^\ast$, the maps
\begin{equation} \label{wstarmod}
  \A \to E^\ast, \quad a \mapsto \left\{
  \begin{array}{c} a \cdot \phi \\ \phi \cdot a \end{array} \right.
\end{equation}
are $w^\ast$-continuous. We write ${\cal Z}^1_{w^\ast}(\A,E^\ast)$ for the 
$w^\ast$-continuous derivations from $\A$ into $E^\ast$. The $w^\ast$-continuity
of the maps (\ref{wstarmod}) implies that ${\cal B}^1(\A,E^\ast) \subset
{\cal Z}^1_{w^\ast}(\A,E^\ast)$, so that ${\cal H}_{w^\ast}^1(\A,E^\ast) :=
{\cal Z}^1_{w^\ast}(\A,E^\ast) / {\cal B}^1(\A,E^\ast)$ is a meaningful definition.
\begin{definition} \label{conam}
A dual Banach algebra $\A$ is {\it Connes-amenable\/} if
${\cal H}^1_{w^\ast}(\A,E^\ast) = \{ 0 \}$ for every $w^\ast$-Banach $\A$-bimodule
$E^\ast$.
\end{definition}
\begin{remarks}
\item Although the notion of Connes-amenability was introduced in \cite{JKR} (for $W^\ast$-algebras), it is most commonly
associated with A.\ Connes' paper \cite{Con76}: This motivates our choice of terminology
(compare also \cite{Hel2}).
\item Definition \ref{conam} is a special case of a notion of amenability introduced in
\cite{CG2}. There, for an arbitrary Banach algebra $\A$ and a submodule $\Phi$ of
$\A^\ast$ satisfying certain properties, $\Phi$-amenability is defined. If
$\A$ is a dual Banach algebra, then $\A_\ast$ satisfies all the requirements for
$\Phi$ in \cite{CG2}, and $\A$ is Connes-amenable if and only if it is $\A_\ast$-amenable
in the sense of \cite{CG2}.
\end{remarks}
\section{The r\^ole of the Radon--Nikod\'ym property}
Let $\mathfrak A$ be a dual Banach algebra, and let ${\mathfrak A}_\ast$
be its predual as in Definition \ref{ddef}. Let ${\mathfrak A}_\ast 
\wTensor {\mathfrak A}_\ast$ be the injective tensor product of 
${\mathfrak A}_\ast$ with itself. Then we have a canonical map from
${\mathfrak A} \Tensor {\mathfrak A}$ into $({\mathfrak A}_\ast 
\wTensor {\mathfrak A}_\ast)^\ast$, which has closed range if
$\mathfrak A$ has the bounded approximation property 
(\cite[16.3, Corollary 2]{DF}).
\par
If $\mathfrak A$ is amenable, a naive approach to show that $\mathfrak A$ 
is super-amenable would be as follows:
\begin{description}
\item[Step 1] Let $( m_\alpha )_\alpha$ be an approximate diagonal for
$\mathfrak A$, and choose an accumulation point $m$ of $( m_\alpha )_\alpha$
in the topology induced by ${\mathfrak A}_\ast \wTensor {\mathfrak A}_\ast$.
\item[Step 2] Show that $m$ is a diagonal for $\mathfrak A$.
\end{description} 
\par
There are problems in both steps (and since there are amenable, dual Banach
algebras which are not super-amenable this is no surprise). In Step 1, the
main problem is that the accumulation point $m \in
({\mathfrak A}_\ast \wTensor {\mathfrak A}_\ast)^\ast$ need not lie
in ${\mathfrak A} \Tensor {\mathfrak A}$. In view of 
\cite[16.6, Theorem]{DF}, it is clear that in order to make Step 1 work,
we have to require the Radon--Nikod\'ym property for $\mathfrak A$.
\par
We have the following theorem: 
\begin{theorem} \label{RNP}
Let ${\mathfrak A}$ be an amenable, dual Banach algebra having both the 
approximation property and the Radon--Nikod\'ym property. Suppose further, 
that there is a family $( I_\lambda )_\lambda$ of $w^\ast$-closed ideals of 
${\mathfrak A}$, each with finite codimension, such that 
$\bigcap_\lambda I_\lambda = \{ 0 \}$. Then
there are $n_1, \ldots, n_k \in \posints$ such that
\[
  {\mathfrak A} \cong \mat_{n_1} \oplus \cdots \oplus \mat_{n_k}.
\]
\end{theorem}
\begin{proof}
Let ${\mathfrak A}_\ast$ denote the predual of ${\mathfrak A}$. Since 
$\mathfrak A$ has both the approximation property and the Radon--Nikod\'ym
property, we have
\[
  {\mathfrak A} \Tensor {\mathfrak A} \cong ({\mathfrak A}_\ast \wTensor
  {\mathfrak A}_\ast )^\ast
\]
by \cite[16.6, Theorem]{DF}. We thus have a natural $w^\ast$-topology on
${\mathfrak A} \Tensor {\mathfrak A}$. Let $( m_\alpha )_\alpha$ be
an approximate diagonal for $\mathfrak A$, and let $m \in {\mathfrak A}
\Tensor {\mathfrak A}$ be a $w^\ast$-accumulation point of
$( m_\alpha )_\alpha$; passing to a subnet we can assume that 
$m = \mbox{$w^\ast$-}\lim_\alpha m_\alpha$.
\par
We claim that $m$ is a diagonal for $\mathfrak A$. It is clear that
$m \in {\cal Z}^0({\mathfrak A}, {\mathfrak A} \Tensor {\mathfrak A})$,
so that all we have to show is that $\Delta m = e_{\mathfrak A}$. Let
$\pi_\lambda \!: {\mathfrak A} \to {\mathfrak A} / I_\lambda$ be the
canonical epimorphism. Since $I_\lambda$ is $w^\ast$-closed, each
quotient algebra ${\mathfrak A} / I_\lambda$ is again dual with the predual
\[
  {^\perp I_\lambda} = \{ \phi \in {\mathfrak A}_\ast :
  \text{$\langle \phi, a \rangle = 0$ for all $a \in I_\lambda$} \}.
\] 
Let $\iota_\lambda \!: {^\perp I_\lambda} \to {\mathfrak A}_\ast$ be
the inclusion map. Then $\pi_\lambda \tensor \pi_\lambda \!:
{\mathfrak A} \Tensor {\mathfrak A} \to {\mathfrak A} / I_\lambda \Tensor 
{\mathfrak A} / I_\lambda$ is the tranpose of $\iota_\lambda \tensor
\iota_\lambda \!: {^\perp I_\lambda} \wTensor {^\perp I_\lambda}
\to {\mathfrak A}_\ast \wTensor {\mathfrak A}_\ast$ (since $I_\lambda$
has finite codimension, we have clearly ${\mathfrak A} / I_\lambda \Tensor 
{\mathfrak A} / I_\lambda\ \cong ({^\perp I_\lambda} \wTensor 
{^\perp I_\lambda})^\ast$). Thus, $\pi_\lambda \tensor \pi_\lambda$ is
$w^\ast$-continuous, so that
\[
  (\pi_\lambda \tensor \pi_\lambda)m = \mbox{$w^\ast$-}\lim_\alpha 
  (\pi_\lambda \tensor \pi_\lambda)m_\alpha.
\]
Since ${\mathfrak A} / I_\lambda \Tensor {\mathfrak A} / I_\lambda$ is
finite-dimensional, there is only one vector space topology on it; in
particular, $(\pi_\lambda \tensor \pi_\lambda)d$ is the norm limit of
$((\pi_\lambda \tensor \pi_\lambda) d_\alpha )_\alpha$. Since
\[
  \Delta_{{\mathfrak A}/I_\lambda} \circ (\pi_\lambda \tensor \pi_\lambda) =
  \pi_\lambda \circ \Delta_{\mathfrak A},
\]
we obtain
\[
  (\pi_\lambda \circ \Delta_{\mathfrak A})m = 
  \lim_\alpha (\Delta_{{\mathfrak A}/I_\lambda} \circ (\pi_\lambda \tensor
  \pi_\lambda)) m_\alpha = e_{{\mathfrak A}/I_\lambda}.
\]
\par
Since $( \pi_\lambda)_\lambda$ separates the points of $\mathfrak A$, it
follows that $\Delta_{\mathfrak A} m = e_{\mathfrak A}$. Hence, $m$ is a
diagonal for $\mathfrak A$.
\end{proof}
\begin{remark}
It is essential for Theorem \ref{RNP} to hold that $\mathfrak A$ has the
Radon--Nikod\'ym property. For example, the algebra
$\ell^\infty$ is an amenable, dual Banach algebra which has a family 
$(I_\lambda)_\lambda$ of $w^\ast$-closed ideals as in Theorem \ref{RNP}, 
but is infinite-dimensional.
\end{remark}
\par
Since separable dual spaces as well as reflexive Banach spaces 
automatically have the Radon--Nikod\'ym property
(\cite[D3]{DF}), we obtain the following corollaries, the first of which 
is a nice dichotomy and the second of which improves
\cite[Corollary 2.3]{GRW} and is related to \cite[Proposition 2.3]{Run}:
\begin{corollary}
Let ${\mathfrak A}$ be an amenable dual Banach algebra having the approximation
property. Suppose further, that there is a family $( I_\lambda )_\lambda$ of 
$w^\ast$-closed ideals of ${\mathfrak A}$, each with finite codimension, such 
that $\bigcap_\lambda I_\lambda = \{ 0 \}$. Then one of the following
holds:
\begin{items}
\item $\mathfrak A$ is not separable;
\item there are $n_1, \ldots, n_k \in \posints$ such that
\[
  {\mathfrak A} \cong \mat_{n_1} \oplus \cdots \oplus \mat_{n_k}.
\]
\end{items}
\end{corollary}
\begin{corollary} \label{refl}
Let ${\mathfrak A}$ be an amenable, reflexive Banach algebra having the 
approximation property. Suppose further, that there is a family 
$( I_\lambda )_\lambda$ of closed ideals of ${\mathfrak A}$, each with finite 
codimension, such that $\bigcap_\lambda I_\lambda = \{ 0 \}$. Then
there are $n_1, \ldots, n_k \in \posints$ such that
\[
  {\mathfrak A} \cong \mat_{n_1} \oplus \cdots \oplus \mat_{n_k}.
\]
\end{corollary}
\begin{remark}
In Corollary \ref{refl} we can replace the hypothesis
that there is a family $( I_\lambda )_\lambda$ of closed ideals of 
${\mathfrak A}$, each with finite codimension, such that 
$\bigcap_\lambda I_\lambda = \{ 0 \}$ by a weaker one. If we assume that
the almost periodic functionals on $\mathfrak A$ separate points, we
still get the same conclusion (this is proved in the same way as
\cite[Proposition 3.1]{Run}). For examples of almost periodic functionals
that do not arise from finite-dimensional quotients, see \cite{DUe}.
\end{remark}
\section{Connes-amenability of biduals}
In this section, we investigate how, for an Arens regular Banach algebra $\A$, the amenability of $\A$ and the Connes-amenability of 
$\A^{\ast\ast}$ are related.
\par
We begin our discussion with some elementary propositions:
\begin{proposition} \label{Cprop1}
Let $\A$ be a Connes-amenable, dual Banach algebra. Then $\A$ has an identity.
\end{proposition}
\begin{proof}
Let $A$ be the Banach $\A$-bimodule whose underlying linear space is $\A$ equipped
with the following module operations:
\[
  a\cdot x := ax \quad\text{and}\quad x\cdot a := 0 \qquad (a,x \in \A).
\]
Obviously, $A$ is a $w^\ast$-Banach $\A$-bimodule the identity on $\A$ into a 
$w^\ast$-continuous derivation. Since ${\cal H}^1_{w^\ast}(\A,A) = \{ 0 \}$,
this means that $\A$ has a right identity. Analoguously, one sees that $\A$ has also a
left identity.
\end{proof}
\begin{proposition} \label{Cprop2}
Let $\A$ be a Banach algebra, let $\B$ be a dual Banach algebra, and let
$\theta \!: \A \to \B$ be a continuous homomorphism with $w^\ast$-dense range. Then:
\begin{items}
\item If $\A$ is amenable, then $\B$ is Connes-amenable.
\item If $\A$ is dual and Connes-amenable, and if $\theta$ is $w^\ast$-continuous,
then $\B$ is Connes-amenable. 
\end{items}
\end{proposition}
\begin{proof}
Immediate from the definitions.
\end{proof}
\begin{corollary} \label{Ccor1}
Let $\A$ be an Arens regular Banach algebra. Then, if $\A$ is amenable, $\A^{\ast\ast}$ is
Connes-amenable.
\end{corollary}
\par
If $\A$ is a $\cstar$-algebra, then the converse of Corollary \ref{Ccor1} holds:
If $\A^{\ast\ast}$ is Connes-amenable, then $\A$ is amenable (\cite{Con76}, \cite{BP}, 
\cite{Haa}, \cite{Eff}, \cite{EK}). This is a deep, specifically $\cstar$-algebraic result, for which no
analogue in the general Banach algebra setting is available (yet). Under certain 
circumstances,
however, a converse of Corollary \ref{Ccor1} holds for general Banach algebras:
\begin{theorem} \label{biduals}
Let $\A$ be an Arens regular Banach algebra which is an ideal in $\A^{\ast\ast}$.
Then the following are equivalent:
\begin{items}
\item $\A$ is amenable.
\item $\A^{\ast\ast}$ is Connes-amenable.
\end{items}
\end{theorem}
\begin{proof}
Since $\A^{\ast\ast}$ is Connes-amenable, it has an identity by Proposition \ref{Cprop1}.
By \cite[Proposition 5.1.8]{Theodore}, this means that $\A$ has a bounded approximate
identity, $(e_\alpha)_\alpha$ say. By \cite{Joh1}, it is therefore sufficient for
$\A$ to be amenable that ${\cal H}^1(\A,E^\ast) = \{ 0 \}$ for each essential
Banach $\A$-bimodule.
\par
Let $E$ be an essential Banach $\A$-bimodule, and let $D \!: A \to E^\ast$ be a
derivation. The following construction is carried out in \cite{Joh1} with the
double centralizer algebra instead of $\A^{\ast\ast}$, but an inspection of the argument
there shows that it carries over to our situation. Since $E$ is essential, there are,
for each $x \in E$, elements $b,c \in \A$ and $y,z \in E$ with $x = b\cdot y = z \cdot c$.
Define an $\A$-bimodule action of $\A^{\ast\ast}$ on $E$, by letting
\[
  a \cdot (b \cdot y) := ab \cdot y \quad\text{and}\quad
  (z \cdot c) \cdot a := z \cdot ca \qquad (a \in \A^{\ast\ast}, \, b,c \in \A, \, y,z \in E).
\]
It can be shown that this module action is well-defined and turns $E$ into a Banach
$\A^{\ast\ast}$-bimodule. Consequently, $E^\ast$ equipped with the corresponding dual
module action is a Banach $\A^{\ast\ast}$-bimodule as well.
\par
We claim that $E^\ast$ is even a $w^\ast$-Banach $\A^{\ast\ast}$-bimodule. Let 
$( a_\beta )_\beta$ be a net in $\A^{\ast\ast}$ such that $a_\beta \stackrel{w^\ast}{\to} 0$,
let $\phi \in E^\ast$, and let $x \in E$. Let $b \in \A$ and $y \in E$ such that $x = 
y \cdot b$. Since the $w^\ast$-topology of $\A^{\ast\ast}$ restricted to $\A$ is the
weak topology, we $b a_\beta \stackrel{w}{\to} 0$, so that
\[
  x \cdot a_\beta = y \cdot b a_\beta \stackrel{w}{\to} 0
\]
and consequently
\[
  \langle x, a_\beta \cdot \phi \rangle = \langle x \cdot a_\beta, \phi \rangle \to 0.
\]
Since $x \in E$ was arbitrary, this means that $a_\beta \cdot \phi \stackrel{w^\ast}{\to} 0$.
Analoguously, one shows that  $\phi\cdot a_\beta  \stackrel{w^\ast}{\to} 0$.
\par
Following again \cite{Joh1}, we define a derivation $\tilde{D} \!: \A^{\ast\ast} \to
E^\ast$ by letting
\[
  \tilde{D} a = \mbox{$w^\ast$-}\lim_\alpha [D(ae_\alpha) - a \cdot De_\alpha].
\]
We claim that $\tilde{D} \in {\cal Z}^1_{w^\ast}(\A^{\ast\ast},E^\ast)$. Let again $(a_\beta)_\beta$
be a net in $\A^{\ast\ast}$ such that $a_\beta \stackrel{w^\ast}{\to} 0$, let $x
\in E$, and let $b \in \A$ and $y \in E$ such that $x = b \cdot y$. Then we have:
\begin{eqnarray*}
  \langle x, \tilde{D}a_\beta \rangle & = & \langle b \cdot y, \tilde{D} a_\beta \rangle \\
  & = & \langle y, (\tilde{D} a_\beta) \cdot b \rangle \\
  & = & \langle y, D(a_\beta b) - a_\beta \cdot Db \rangle \\
  & \to & 0, 
\end{eqnarray*}
since $D$ is weakly continuous and $E^\ast$ is a $w^\ast$-Banach $\A^{\ast\ast}$-bimodule.
\par
From the Connes-amenability of $\A^{\ast\ast}$ we conclude that $\tilde{D}$, and hence
$D$, is inner.
\end{proof}
\begin{remark}
In \cite{Gou}, F.\ Gourdeau showed that, whenever $\A$ is a Banach algebra, $E$ is a Banach
$\A$-bimodule, and $D \!: \A \to E$ is a derivation, there is an $\A^{\ast\ast}$-bimodule
action on $E^{\ast\ast}$, turning $D^{\ast\ast} \!: \A^{\ast\ast} \to E^{\ast\ast}$ into
a (necessarily $w^\ast$-continuous) derivation. However, even if $E$ is a dual 
Banach $\A$-bimodule, there is no need for $E^{\ast\ast}$ to be a $w^\ast$-Banach
$\A$-bimodule, so that, in general, we cannot draw any conclusion on the amenability of $\A$
from the Connes-amenability of $\A^{\ast\ast}$.
\end{remark}
\par
By \cite[Theorem 6.9]{GJW}, the space $\ell^p \oplus \ell^q$ with
$p,q \in (1,\infty) \setminus \{ 2 \}$ and $p \neq q$ has the property that
${\cal K}(\ell^p \oplus \ell^q)$ is not amenable. Hence, Theorem \ref{biduals} 
yields:
\begin{corollary}
Let $p,q \in (1,\infty) \setminus \{ 2 \}$, $p \neq q$. Then ${\cal L}(\ell^p \oplus 
\ell^q)$ is not Connes-amenable.
\end{corollary}
\begin{proof}
Since ${\cal K}(\ell^p \oplus \ell^q)^{\ast\ast} \cong {\cal L}(\ell^p \oplus \ell^q)$, and
since ${\cal K}(\ell^p \oplus \ell^q)$ is not amenable, ${\cal L}(\ell^p \oplus \ell^q)$
is not Connes-amenable by Theorem \ref{biduals}.
\end{proof}
\par
Let $\A$ be a dual Banach algebra, and let $E$ be a Banach $\A$-bimodule. Then we
call an element $\phi \in E^\ast$ a {\it $w^\ast$-element\/} if the maps
(\ref{wstarmod}) are $w^\ast$-continuous.
\begin{definition} \label{sconam}
A dual Banach algebra with identity $\A$ is called {\it strongly Connes-amenable\/} if, for each
each unital Banach $\A$-bimodule $E$, every $w^\ast$-continuous derivation $D \!:
\A \to E^\ast$ whose range consists
of $w^\ast$-elements is inner.
\end{definition}
\par
We shall give an instrinsic characterization of strongly Connes-amenable dual Banach
algebras similar to the one given in \cite{Joh2} for amenable Banach algebras.
Recall a few definitions from \cite{CG2} (with a different notation, however). 
Let $\A$ be a dual Banach algebra with identity, and let ${\cal L}^2_{w^\ast}(\A,\comps)$ be the 
space of separately 
$w^\ast$-continuous bilinear functionals on $\A$. Clearly, ${\cal L}^2_{w^\ast}(\A,\comps)$ is 
Banach $\A$-submodule of ${\cal L}^2(\A,\comps) \cong (\A \ptensor \A)^\ast$. Define
\[
  (\A \ptensor_{w^\ast} \A)^{\ast\ast} := {\cal L}^2_{w^\ast}(\A,\comps)^\ast.
\]
Note that the notation $(\A \ptensor_{w^\ast} \A)^{\ast\ast}$ is merely symbolic: In
general, $(\A \ptensor_{w^\ast} \A)^{\ast\ast}$ is not a bidual space. There is a
canonical embedding of the algebraic tensor product $\A \tensor \A$ into 
$(\A \ptensor_{w^\ast} \A)^{\ast\ast}$, so that we may identify $\A \tensor \A$ with
a submodule of $(\A \ptensor_{w^\ast} \A)^{\ast\ast}$. It is immediate that 
$\A \tensor \A$ consists of $w^\ast$-elements of $(\A \ptensor_{w^\ast} \A)^{\ast\ast}$.
Since multiplicaton in a dual Banach algebra is separately $w^\ast$-continuous, we have
\[
  \Delta^\ast \A_\ast \subset {\cal L}^2_{w^\ast}(\A,\comps),
\]
so that the multiplication operator $\Delta$ on $\A \tensor \A$ extends to 
$(\A \ptensor_{w^\ast} \A)^{\ast\ast}$; we shall denote this extension by 
$\Delta_{w^\ast}^{\ast\ast}$. A {\it virtual $w^\ast$-diagonal\/} for $\A$ (in the
terminology of \cite{CG2}: an $\A_\ast$-virtual diagonal) is an element $M \in
(\A \ptensor_{w^\ast} \A)^{\ast\ast}$ such that
\[
  a \cdot M = M \cdot a   
  \quad (a \in \A) \qquad\text{and}\qquad \Delta^{\ast\ast}_{w^\ast} M = e_\A.
\]
\par
In \cite{CG2}, G.\ Corach and J.\ E.\ Gal\'e showed that a dual Banach algebra with
a virtual $w^\ast$-diagonal is necessarily Connes-amenable, and wondered if the converse
was also true. For strong Connes-amenability, the corresponding question is 
easy to answer:
\begin{theorem} \label{coneq}
For a dual Banach algebra $\A$, the following are equivalent:
\begin{items}
\item $\A$ has a virtual $w^\ast$-diagonal.
\item $\A$ is strongly Connes-amenable.
\end{items}
\end{theorem}
\begin{proof}
On \cite[p.\ 90]{CG2}, it is shown that (i) implies the Connes-amenability of $\A$ (the
argument for von Neumann algebras from \cite{Eff} carries over verbatim). A
closer inspection of the argument given there, however, shows that we already obtain
strong Connes-amenability.
\par
For the converse, consider the derivation $\ad_{e_\A \tensor e_\A}$. Then, clearly, 
$\ad_{e_\A \tensor e_\A}$
attains its values in the $w^\ast$-elements of $\ker \Delta_{w^\ast}^{\ast\ast}$.
By the definition of strong Connes-amenability, there is $N \in 
\ker \Delta_{w^\ast}^{\ast\ast}$ such that $\ad_N = \ad_{e_\A \tensor e_\A}$. It follows
that $D := e_\A \tensor e_\A - N$ is a virtual $w^\ast$-dagonal for $\A$.
\end{proof}
\begin{remark}
In \cite{Eff}, E.\ G.\ Effros proves that a von Neumann algebra is Connes-amenable if and only
if it has a virtual $w^\ast$-diagonal. Hence, von Neumann algebras are Connes-amenable if and only
if they are strongly Connes-amenable.
\end{remark}
\par
For certain Banach algebras $\A$, the strong Connes-amenability of $\A^{\ast\ast}$ entails the
amenability of $\A$:
\begin{theorem} \label{coneq2}
Let $\A$ be a Banach algebra with the following properties:
\begin{items}
\item Every bounded linear map from $\A$ to $\A^\ast$ is weakly compact.
\item $\A^{\ast\ast}$ is strongly Connes-amenable.
\end{items}
Then $\A$ is amenable.
\end{theorem}
\begin{proof}
By \cite[1.29 Satz]{Gro}, condition (i) implies (and is, in fact, equivalent to) that
every bounded, bilinear map from $\A \times \A$ into any Banach space is Arens regular;
in particular, it ensures that $\A^{\ast\ast}$ is indeed a dual Banach algebra. It is thus
an immediate consequence of (i) that 
\begin{equation} \label{modiso}
  (\A \Tensor \A)^{\ast\ast} \cong (\A^{\ast\ast} \Tensor_{w^\ast} \A^{\ast\ast})^{\ast\ast},
\end{equation}
as Banach $\A$-bimodules. Since $\A^{\ast\ast}$ has a virtual $w^\ast$-diagonal by Theorem
\ref{coneq}, the isomorphism (\ref{modiso}) ensures the existence of a virtual diagonal for
$\A$. Thus, $\A$ is amenable.
\end{proof}
\begin{examples}
\item Every $\cstar$-algebra satisfies Theorem \ref{coneq2}(i).
\item Let $E$ be a reflexive Banach space with an unconditional basis. It is implicitly proved
in \cite{Ulg} (although not explicitly stated) that ${\cal K}(E)$ satisfies Theorem 
\ref{coneq2}(i).
\end{examples}
\section{Dual Banach algebras associated with locally compact groups}
As was shown by A.\ Connes (\cite{Con}), Connes-amenable von Neumann algebras (with
separable predual) are injective. An alternative proof for this is given in \cite{BP}. We
now give an analogue --- in the spirit of \cite[Proposition 2.2]{CG2} and
\cite[Proposition 2.2]{CG1} --- for arbitrary Connes-amenable, dual Banach
algebras.
\par
If $S$ is any subset of an algebra $\B$, we use $Z_\B(S)$ to denote the centralizer
of $S$ in $\B$, i.e.\
\[
  Z_\B(S) := \{ b \in \B : \text{$bs = sb$ for all $s \in S$} \}.
\]
In case $\B = {\cal L}(E)$ for some Banach space $E$, we also write $S'$ instead of $Z_\B(S)$. 
Recall (from \cite{BP} or \cite{CG1}) the definition of a {\it quasi-expectation\/}: If
$\B$ is a Banach algebra, and if $\A$ is a closed subalgebra of $\B$, a
quasi-expectation is a bounded projection ${\cal Q} \!: \B \to \A$ satisfying
\[
  {\cal Q}(axb) = a ({\cal Q}x)b \qquad (a,b \in \A, \, x \in \B).
\]
Note that we do not require $\| {\cal Q} \| = 1$.
\begin{theorem} \label{quasexp}
Let $\A$ be a Connes-amenable dual Banach algebra, let $\B$ be a dual Banach algebra with
identity, and let $\theta \!: \A \to \B$ be a unital, $w^\ast$-continuous homomorphism.
Then there is a quasi-expectation ${\cal Q}\!: \B \to Z_{\B}(\theta(\A))$.
\end{theorem}
\begin{proof}
More or less a verbatim copy of the proof of \cite[Theorem 3]{BP}.    
\end{proof}    
\par
In this section, we will use Theorem \ref{quasexp} to
characterize the Connes-amenability of some dual Banach algebra which arise naturally
in abstract harmonic analysis.
\par
For non-discrete, abelian $G$, it has long been known that there are non-zero
point derivations on $M(G)$ (\cite{BM}), so that $M(G)$ cannot be amenable.
In \cite{LL}, A.\ T.-M.\ Lau and R.\ J.\ Loy started investigating the amenability of
$M(G)$ for certain non-abelian $G$. In particular, they were able to show that,
for connected $G$, the algebra $M(G)$ is amenable only if $G = \{ e \}$
(\cite[Theorem 2.4]{LL}). Ultimately, H.\ G.\ Dales, F.\ Gharamani, and A.\ Ya.\
Helemski\u{\i} (\cite{DGH}) proved: The measure algebra $M(G)$ is amenable if and only if
$G$ is discrete --- so that $M(G) = \ell^1(G) = L^1(G)$ --- and amenable. 
\par
The picture is completely different for Connes-amenability. We have, for example:
\begin{proposition} \label{compact}
Let $G$ be a compact group. Then $M(G)$ is strongly Connes-amenable.
\end{proposition}
\begin{proof}
By Theorem \ref{coneq}, it is sufficient to construct a virtual $w^\ast$-diagonal for $M(G)$.
\par
For $\phi \in \in {\cal L}_{w^\ast}^2(M(G),\comps)$, define $\bar{\phi} \!: G \times G \to \comps$
and $\tilde{\phi} \!: G \to \comps$ through
\[
  \bar{\phi}(x,y) := \phi(\delta_x, \delta_y ) \quad\text{and}\quad
  \tilde{\phi}(x) := \bar{\phi}(x,x^{-1}) \qquad (x,y \in G).
\]
Then $\bar{\phi}$ is separately continuous on $G \times G$ and thus belongs to 
$L^\infty(G \times G,\mu \times \nu)$ for any $\nu, \mu \in M(G)$ (\cite{Joh0}). Since 
(normalized) Haar measure belongs to $M(G)$, this implies that $\tilde{\phi} \in L^\infty(G) 
\subset L^1(G)$.
\par
Let $m$ denote normalized Haar measure on $G$, and define $M \in
(M(G) \ptensor_{w^\ast} M(G))^{\ast\ast}$ via
\[
  \langle \phi, M \rangle := \int_G \tilde{\phi}(x) \, dm(x)
  \qquad (\phi \in {\cal L}_{w^\ast}^2(M(G),\comps)).
\]
It is routinely checked that $\Delta^{\ast\ast}_{w^\ast} M = \delta_e$. Let
$\mu \in M(G)$, and let $\phi \in {\cal L}^2_{w^\ast}(M(G),\comps)$. Then we have:
\begin{eqnarray} 
  \langle \phi, \mu \cdot M \rangle & = & \langle \phi \cdot \mu , M \rangle \nonumber \\
  & = & \int_G \widetilde{\phi \cdot \mu}(x) \, dm(x) \nonumber \\
  & = & \int_G \overline{\phi \cdot \mu}(x,x^{-1}) \, dm(x) \nonumber \\
  & = & \int_G \int_G \bar{\phi}(yx,x^{-1}) \, d\mu(y) \, dm(x) \nonumber \\ 
  & = & \int_G \int_G \bar{\phi}(yx,x^{-1}) \, dm(x) \, d\mu(y), 
  \qquad\text{by Fubini's theorem}, \label{Fubini1} \\
  & = & \int_G \int_G \bar{\phi}(x,x^{-1}y) \, dm(x) \, d\mu(y), 
  \qquad\text{substituting $y^{-1}x$ for $x$}, \nonumber \\
  & = & \int_G \int_G \bar{\phi}(x,x^{-1}y ) \, d\mu(y) \, dm(x), 
  \qquad\text{again by Fubini's theorem}, \label{Fubini2} \\
  & = & \int_G \widetilde{\mu \cdot \phi}(x) \, dm(x) \nonumber \\
  & = & \langle \mu \cdot \phi, M \rangle \nonumber \\
  & = & \langle \phi, M \cdot \mu \rangle. \nonumber
\end{eqnarray}
Thus, $M$ is virtual $w^\ast$-diagonal for $M(G)$.
\end{proof}
\begin{remarks}
\item Most of the proof of Proposition \ref{compact} still works for not necessarily compact, amenable
$G$, where we replace Haar measure in the definiton of $M$ by a left invariant mean on
$L^\infty(G)$. In this, more general situation, however, we have no substitute for the two
applications (\ref{Fubini1}) and (\ref{Fubini2}) of Fubini's theorem.
\item The locally compact groups $G$ for which $M(G)$ is Connes-amenable are characterized in
the companion paper \cite{Run2}: They are precisely the amenable groups. In the same paper, it is also
shown that $M(G)$ has a virtual $w^\ast$-diagonal if and only if it is Connes-amenable.
\end{remarks}
\par
Another important dual Banach algebra associated with a locally compact group $G$ is
the group von Neumann algebra $VN(G)$. Its Connes-amenability was characterized in
terms of $G$ by A.\ T.-M.\ Lau and A.\ L.\ T.\ Paterson (\cite{LP}). There are analogues
of $VN(G)$ acting on $L^p(G)$ for $p \in (1,\infty)$, the so-called algebras of
{\it pseudo-mesures\/} $PM_p(G)$ (for information on these algebras and further 
references, see \cite{Eym}). We shall now prove an extension of \cite[Corollary 3.3]{LP} for 
these algebras of pseudo-measures.
\par
Recall (from \cite{LP}, for example) that a locally compact group is {\it inner
amenable\/} if there is a state $m$ on $L^\infty(G)$ such that
\[
  \langle \delta_x \ast \phi \ast \delta_{x^{-1}}, m \rangle = \langle \phi, m \rangle
  \qquad (x \in G, \, \phi \in L^\infty(G)).
\]
Every amenable group is inner amenable, but so is evey $[IN]$-group.
\begin{theorem} \label{pseudothm}
For a locally compact group $G$ consider the following:
\begin{items}
\item $G$ is amenable.
\item $M(G)$ is Connes-amenable.
\item $PM_p(G)$ is Connes-amenable for every $p \in (1,\infty)$.
\item $VN(G)$ is Connes-amenable.
\item $PM_p(G)$ is Connes-amenable for one $p \in (1,\infty)$.
\end{items}
Then we have:
\begin{center}
  {\rm (i)} $\Longrightarrow$  {\rm (ii)} $\Longrightarrow$  {\rm (iii)} $\Longrightarrow$  
  {\rm (iv)} $\Longrightarrow$  {\rm (v)}.
\end{center}
If $G$ is inner amenable,\ {\rm (v)} $\Longrightarrow$ {\rm (i)} holds, too.
\end{theorem}
\begin{proof}
(i) $\Longrightarrow$ (ii): This is is clear in view of Proposition \ref{Cprop2}(i).
\par
(ii) $\Longrightarrow$ (iii): This follows from Proposition \ref{Cprop2}(ii).
\par
(iii) $\Longrightarrow$ (iv) $\Longrightarrow$ (v): Since $VN(G) = PM_2(G)$,
this is obvious.
\par
(v) $\Longrightarrow$ (i) for $G$ inner amenable: 
For any $r \in [1,\infty)$, let $\lambda_r$ and $\rho_r$
denote the regular left and right representation, respectively, of $G$ on $L^r(G)$.
\par
By \cite[Proposition 1]{LR}, it follows from the inner amenability 
of $G$ that there is a net $(f_\alpha)_\alpha$ of positive $L^1$-functions with
$\| f_\alpha \|_1 = 1$ such that 
\[
  \| \delta_x \ast f_\alpha \ast \delta_{x^{-1}} - f_\alpha \|_1 \to 0 \qquad
  (x \in G), 
\]
or equivalently
\begin{equation} \label{innam}
   \| \lambda_1(x^{-1}) f_\alpha - \rho_1(x)_1f_\alpha \|_1 \to 0 \qquad (x \in G).
\end{equation} 
Let $q \in (1,\infty)$ be the index dual to $p$. Let $\xi_\alpha := f_\alpha^{1/p}$, and
let $\eta_\alpha := f_\alpha^{1/q}$, so that $\xi_\alpha \in L^p(G)$ and $\eta_\alpha
\in L^q(G)$. It follows from (\ref{innam}) and \cite[4.3(1)]{Pat} that
\begin{equation} \label{inneram}  
  \left. \begin{array}{c}
  \| \lambda_p(x^{-1}) \xi_\alpha - \rho_p(x)\xi_\alpha \|_p \to 0 \\
   \text{and} \\
  \| \lambda_q(x^{-1}) \eta_\alpha - \rho_q(x)\eta_\alpha \|_q \to 0
  \end{array} \right\} \qquad (x \in G).
\end{equation}
For $\phi \in UC(G)$, let $M_\phi \in {\cal L}(L^p(G))$
be defined by pointwise multiplication with $\phi$. By Theorem \ref{quasexp}
--- applied to $\A = PM_p(G)$, $\B = {\cal L}(L^p(G))$, and $\theta$ the canonical 
representation of $PM_p(G)$ on $L^p(G)$ --- , there is a quasi-expectation 
${\cal Q} \!: {\cal L}(L^p(G)) \to PM_p(G)'$. Define $m_\alpha \in UC(G)^\ast$ by letting
\[
  \langle \phi, m_\alpha \rangle := 
  \langle {\cal Q}(M_\phi)\xi_\alpha,\eta_\alpha \rangle
  \qquad (\phi \in UC(G)).
\] 
Let $\cal U$ be an ultrafilter on the index set of $( f_\alpha)_\alpha$ that dominates the order
filter, and define
\[
  \langle \phi, m \rangle := \lim_{\cal U} \langle \phi, m_\alpha \rangle
  \qquad (\phi \in UC(G)).
\]
Note that $\rho_p(G) \subset PM_p(G)'$, and observe again that
\[
  \rho_p(x^{-1}) M_\phi \rho_p(x) = M_{\phi \ast \delta_x}
  \qquad (x \in G, \, \phi \in UC(G)).
\]
We then obtain for $x \in G$ and $\phi \in UC(G)$:
\begin{eqnarray*}
  \langle \phi \ast \delta_x, m \rangle & = & 
  \lim_{\cal U} \langle \phi \ast \delta_x, m_\alpha \rangle \\
  & = & \lim_{\cal U} \langle {\cal Q}(\rho_p(x^{-1}) M_\phi \rho_p(x)) 
                        \xi_\alpha, \eta_\alpha \rangle \\
  & = & \lim_{\cal U} \langle \rho_p(x^{-1}) ({\cal Q} M_\phi) \rho_p(x)
                        \xi_\alpha, \eta_\alpha \rangle \\
  & = & \lim_{\cal U} \langle ({\cal Q} M_\phi) \rho_p(x)
                        \xi_\alpha, \rho_p(x) \eta_\alpha \rangle \\ 
  & = & \lim_{\cal U} \langle ({\cal Q} M_\phi) \lambda_p(x^{-1}) \xi_\alpha, 
                              \lambda_p(x^{-1}) \eta_\alpha \rangle,
                      \qquad\text{by (\ref{inneram})}, \\
  & = & \lim_{\cal U} \langle \lambda_p(x) {\cal Q} M_\phi) \lambda_p(x^{-1}) \xi_\alpha, 
                              \eta_\alpha \rangle \\
  & = & \lim_{\cal U} \langle ({\cal Q}M_\phi)\xi_\alpha, \eta_\alpha \rangle \\
  & = & \langle \phi, m \rangle.
\end{eqnarray*}
Hence, $m$ is right invariant. Clearly, $\langle 1,m \rangle = 1$. Taking the
positive part of $m$ and normalizing it, we obtain a right invariant mean on $UC(G)$. 
\end{proof}
\begin{remarks}
\item The hypothesis that $G$ be inner amenable cannot be dropped: As pointed out
on \cite[p.\ 84]{Pat}, $SL(2,\reals)$ is not amenable, but of type I, so that 
$VN(SL(2,\reals))$ is Connes-amenable.
\item In \cite{LP}, Lau and Paterson show that an inner amenable, locally compact group $G$ such
that $VN(G)$ has Schwartz' property $(P)$ is already amenable. Schwartz' property
$(P)$ and Connes-amenability are equivalent (see \cite[(2.35)]{Pat}), but the 
implication from Connes-amenability to Schwartz' property $(P)$ is a deep result
by itself. The proof of Theorem \ref{pseudothm} is thus much simpler, even in the
particular case $p = 2$.
\end{remarks}
\section{A nuclear-free characterization of amenable $W^\ast$-algebras}
As we have mentioned in the introduction, the following characterization 
of the amenable $W^\ast$-algebras is known:
\begin{theorem} \label{wstar}
For a $W^\ast$-algebra $\mathfrak A$, the following are equivalent:
\begin{items}
\item $\mathfrak A$ is amenable.
\item There are  hyperstonean, compact spaces $\Omega_1, \ldots,
\Omega_n$ and $n_1, \ldots, n_k \in \posints$ such that
\[
  {\mathfrak A} \cong 
  \bigoplus_{j=1}^k {\mathbb M}_{n_j} \tensor {\cal C}(\Omega_j).
\]
\end{items}
\end{theorem}
\par
The implication (ii) $\Longrightarrow$ (i) is obvious, and the converse
is a consequence of \cite{Con}, \cite[(1.9) Corollary]{Was}, and the
structure theory of $W^\ast$-algebras.
\par
In this section, we give a proof of Theorem \ref{wstar} which avoids the
amenability-nuclearity nexus and only relies on standard facts about
$W^\ast$-algebras, where we define a standard fact as one that can be
found in one of the standard books on the subject such as
\cite{Dix}, \cite{KR}, \cite{Sak}, and \cite{Tak}. 
\par
In analogy with \cite{Was}, we define:
\begin{definition}
Let $\mathfrak A$ be a Banach $^\ast$-algebra. We say that $\mathfrak A$ 
is of {\it type\/ $(QE)$\/} if, for each $^\ast$-representation 
$(\pi, {\mathfrak H})$, there is a quasi-expectation ${\cal Q} \!: 
{\cal B}({\mathfrak H}) \to \pi({\mathfrak A})^{\prime\prime}$.
\end{definition}
\begin{remark}
It follows from \cite{BP} that every $\cstar$-algebra which is of
type $(QE)$ is already of type $(E)$ in the sense of \cite{Was}. However, 
since we strive to keep the $W^\ast$-theory required for the proof of
Theorem \ref{wstar} to a minimum, we shall not use this fact.
\end{remark}
\par
The next two results are essentially already contained in \cite{Tom}. We
give proofs here requiring a minimum of $W^\ast$-theory:
\begin{lemma} \label{qext}
Let $\mathfrak A$ be a von Neumann algebra acting on a Hilbert space
$\mathfrak H$. Then the following are equivalent:
\begin{items}
\item There is a quasi-expectation ${\cal Q} \!: {\cal B}({\mathfrak H}) \to
{\mathfrak A}$.
\item For every faithful, normal representation $(\pi,{\mathfrak K})$ of
$\mathfrak A$, there is a quasi-expectation ${\cal Q} \!: {\cal B}
({\mathfrak K}) \to \pi({\mathfrak A})^{\prime\prime}$.
\end{items}
\end{lemma}
\begin{proof}
Let $\pi$ be a faithful $W^\ast$-representation of $\M$ on a Hilbert space ${\mathfrak K}$, and let $\N := \pi({\mathfrak \M})$. Using the idea of the proof of 
\cite[Th\'eor\`eme 3, \S I. 4]{Dix}, we can choose a third Hilbert space $\mathfrak L$ such that 
${\mathfrak \M} \tensor \id_{\mathfrak L}$, i.e.\ the algebra $\{ x \tensor \id_{\mathfrak L} : x \in \M \}$, 
on ${\mathfrak H} \bar{\tensor} {\mathfrak L}$ and $\N \tensor \id_{\mathfrak L}$ on ${\mathfrak K} \bar{\tensor} {\mathfrak L}$ are spatially
isomorphic. Fix $\xi_0, \eta_0 \in \mathfrak L$ such that $\langle \xi_0, \eta_0 \rangle = 1$, and define 
${\cal P}_0 \!: {\cal L}({\mathfrak H} \bar{\tensor} {\mathfrak L}) \to {\cal L} ({\mathfrak H})$ through  
\begin{equation} \label{qexp}
   \langle ({\cal P}_0 T) \xi, \eta \rangle := \langle T(\xi \tensor \xi_0), \eta \tensor \eta_0 \rangle
   \qquad (T \in {\cal L}({\mathfrak H}), \, \xi, \eta \in {\mathfrak H}).
\end{equation}
Identifying ${\cal L}({\mathfrak H})$ with ${\cal L}({\mathfrak H}) \tensor \id_{\mathfrak L}$, we see that ${\cal P}_0$ is
a projection onto ${\cal L}({\mathfrak H})$. Furthermore, for $T \in {\cal L}({\mathfrak H} \bar{\tensor} {\mathfrak L})$ and
$R,S \in {\cal L}({\mathfrak H})$, we have
\begin{eqnarray*}
  {\cal P}_0((R \tensor \id_{\mathfrak L})T(S \tensor \id_{\mathfrak L}))\xi, \eta \rangle
  & = & \langle (R \tensor \id_{\mathfrak L})T(S \tensor \id_{\mathfrak L})(\xi \tensor \xi_0), \eta \tensor \eta_0 \rangle \\
  & = & \langle T(S\xi \tensor \xi_0), R^\ast \eta \tensor \eta_0 \rangle \\
  & = & \langle ({\cal P}_0T)S\xi, R^\ast \eta \rangle \\
  & = & \langle R({\cal P}_0T)S\xi, \eta \rangle
  \qquad (\xi, \eta \in {\mathfrak H}),
\end{eqnarray*}
so that ${\cal P}_0$ is a quasi-expectation. Let ${\cal P} \!: {\cal L}({\mathfrak H}) \to \M$ be a quasi-expectation, and define
\[
  \tilde{\cal P} \!: {\cal L}({\mathfrak H} \bar{\tensor} {\mathfrak L}) \to \M \tensor \id_{\mathfrak L},
  \quad T \mapsto ({\cal P} \circ {\cal P}_0) T \tensor \id_{\mathfrak L};
\]
it is clear that $\tilde{\cal P}$ is also a quasi-expectation. Since $\M \tensor \id_{\mathfrak L}$ and $\N \tensor \id_{\mathfrak L}$
are spatially isomorphic, there is a quasi-expecation $\tilde{\cal Q} \!: {\cal L}({\mathfrak K} \bar{\tensor} {\mathfrak L})
\to \N \tensor \id_{\mathfrak L}$. Fix again $\xi_0, \eta_0 \in \mathfrak L$ such that $\langle \xi_0, \eta_0 \rangle = 1$,
and define ${\cal Q}_0 \!: {\cal L}({\mathfrak K} \bar{\tensor} {\mathfrak L}) \to {\cal L}({\mathfrak K})$ as in (\ref{qexp}).
Then
\[
  {\cal Q} \!: {\cal L}({\mathfrak K}) \to \N, \quad T \mapsto ({\cal Q}_0 \circ
 \tilde{\cal Q})(T \tensor \id_{\mathfrak L})
\]
is the desired quasi-expectation.
\end{proof}
\begin{proposition} \label{QE}
Every amenable Banach $^\ast$-algebra is of type $(QE)$.
\end{proposition}
\begin{proof}
Let $(\pi, {\mathfrak H})$ be a $^\ast$-representation of $\mathfrak A$. 
Then, by \cite[Exercise 7.6.46]{KR}, there is faithful, normal, semifinite
weight on the von Neumann algebra $\pi({\mathfrak A})^{\prime\prime}$. Let
$(\rho,{\mathfrak K})$ be the faithful, normal representation of
$\pi({\mathfrak A})^{\prime\prime}$ constructed from this weight 
(\cite[Theorem 7.5.3]{KR}). By \cite[Proposition 2.2]{CG1} --- applied to
the representation $(\rho \circ \pi, {\mathfrak K})$ of $\mathfrak A$ ---,
there is a quasi-expectation ${\cal P} \!: {\cal B}({\mathfrak K}) \to
(\rho \circ \pi)({\mathfrak A})^\prime$. Let $J \!: {\mathfrak K} \to 
{\mathfrak K}$ be the conjugate linear isometry from \cite[Theorem 9.2.37]{KR},
i.e.\ 
\[
  J^2 = \id_{\mathfrak K} \qquad\text{and}\qquad
  J (\rho \circ \pi)({\mathfrak A})' J = 
  (\rho \circ \pi)({\mathfrak A})^{\prime\prime}.
\]
Define
\[
  \tilde{\cal Q} \!: {\cal B}({\mathfrak K}) \to 
  (\rho \circ \pi)({\mathfrak A})^{\prime\prime}, \quad
  a \mapsto J {\cal P}(JaJ) J.
\]
It is immediate that $\tilde{\cal Q}$ is a quasi-expectation. Since $\rho$ is
normal and faithful, it follows from Lemma \ref{qext} that there is a
quasi-expectation ${\cal Q} \!: {\cal B}({\mathfrak H}) \to
\pi({\mathfrak A})^{\prime\prime}$.
\end{proof}
\par
Our next lemma is an analogue of \cite[(1.2) Proposition]{Was}, whose proof
carries over almost verbatim:
\begin{lemma} \label{QElemma}
Let $\mathfrak A$ be a $\cstar$-algebra of type $(QE)$, and let 
$\mathfrak B$ be a $\cstar$-algebra such that there is a quasi-expectation
${\cal Q} \!: {\mathfrak A} \to {\mathfrak B}$. Then $\mathfrak B$ is of
type $(QE)$.
\end{lemma}
\par
Also, the following proposition is well known: It is an immediate consequence of
\cite[Corollary 3.3]{LP}, \cite[Theorem 7.2]{Tom}, amd \cite[Theorem 2]{BP}.
We prefer, however, to indicate a proof that completely avoids $W^\ast$-theoretical
arguments and uses ideas from \cite{CG1} instead:
\begin{proposition} \label{amgroups}
For an inner amenable group $G$, the following are equivalent:
\begin{items}
\item $G$ is amenable.
\item There is a quasi-expectation ${\cal Q} \!: {\cal B}(L^2(G)) \to
VN(G)$.
\end{items}
\end{proposition}
\begin{proof}
For the proof of (i) $\Longrightarrow$ (ii) note that $VN(G) = 
\rho_2(G)^\prime$. An application of \cite[Proposition 2.2]{CG1} then yields the claim.
\par
For the converse implication, we proceed as in the proof of
Theorem \ref{pseudothm}, with the r\^oles of $\lambda_2$ and $\rho_2$
interchanged, and obtain a left-invariant $m \in UC(G)^\ast$ with
$\langle 1, m \rangle =1$. 
\end{proof}
\par
As in \cite{Was}, we obtain:
\begin{corollary}
The $W^\ast$-algebras $VN(\free_2)$ and ${\mathbb M}_\infty$ are not of type
$(QE)$ and thus, in particular, are not amenable.
\end{corollary} 
\begin{proof}
Since $\free_2$ is well known (and easily seen) to be non-amenable
(\cite[(0.6) Example]{Pat}), there is no quasi-expectation 
${\cal Q} \!: {\cal B}(\ell^2(\free_2)) \to VN(\free_2)$ by Proposition
\ref{amgroups}. 
\par
The case of ${\mathbb M}_\infty$ is reduced to $VN(\free_2)$ as in
\cite{Was}: We can find a maximal ideal $M$ of ${\mathbb M}_\infty$ 
corresponding to a point in $\beta\posints \setminus \posints$ such that
${\mathfrak A} := {\mathbb M}_\infty / M$ contains a $W^\ast$-subalgebra
$\mathfrak B$ which is isomorphic to $VN(\free_2)$. As pointed out
in \cite[pp.\ 358--359]{Tak}, $\mathfrak A$ is a type ${\rm II}_1$ factor,
which by \cite[4.4.23 Proposition]{Sak} means that there is a (norm one)
quasi-expectation ${\cal Q} \!: {\mathfrak A} \to {\mathfrak B}$. Hence,
if ${\mathbb M}_\infty$ were of type $(QE)$, the same would be true for
$\mathfrak A$ and, by Lemma \ref{QElemma}, for ${\mathfrak B} \cong 
VN(\free_2)$. However, as we have just seen, $VN(\free_2)$ fails to be of
type $(QE)$.
\end{proof}
\begin{remark}
Although the proof for the non-amenability of ${\mathbb M}_\infty$ given
here is more elementary than the one from \cite{LLW}, it is not yet as
satisfying as we would like it to be: The proof of Proposition \ref{QE}
relies on Tomita--Takesaki theory, which still seems to be far too deep
in order to prove a result on an algebra as plain as ${\mathbb M}_\infty$.
We could have avoided the use of Tomita--Takesaki theory in the proof
of Proposition \ref{QE} by defining type $(QE)$ via 
$\pi({\mathfrak A})^\prime$ instead of $\pi({\mathfrak A})^{\prime\prime}$:
Proposition \ref{QE} would then have been a straightforward application of
\cite[Proposition 2.2]{CG1}. Then, however, we would have required
Tomita--Takesaki theory in the proof of Lemma \ref{QElemma}.
\end{remark}
\pf{Proof of Theorem \ref{wstar}}
Suppose that $\mathfrak A$ is not of the form in Theorem \ref{wstar}(ii).
Then, just as in \cite{Was}, it follows from the general structure theory of
$W^\ast$-algebras that $\mathfrak A$ contains ${\mathbb M}_\infty$ as
a closed subalgebra. As in the proof of \cite[(1.8) Corollary]{Was}, 
a norm one projection ${\cal Q} \!: {\mathfrak A} \to {\mathbb M}_\infty$
can be constructed, which is easily, i.e.\ directly and without 
\cite[Theorem 3.1]{Tom}, seen to be a quasi-expectation. 
\qed
\par
The proof Theorem \ref{wstar} would be considerably easier if we could
show by elementary means that that, for an amenable $\cstar$-algebra
$\mathfrak A$, a $\cstar$-subalgebra $\mathfrak B$, and a quasi-expectation
${\cal Q} \!: {\mathfrak A} \to {\mathfrak B}$, the algebra $\mathfrak B$ is
amenable as well. This is indeed true, but in order to prove it we need the
deep connections between amenability, nuclearity and injectivity. 
\par
The corresponding claim for general Banach algebras is false:
\begin{example}
Let $G$ be a compact group. Then $L^1(G)$ is amenable, so that $L^1(G)$ has
an approximate diagonal $( m_\alpha )_\alpha$. Let $\cal U$ be an
ultrafilter on the index set of $(d_\alpha )_\alpha$ which dominates the order
filter, and define 
\[
  {\cal Q} \!: L^1(G) \to Z(L^1(G)), \quad f \mapsto \lim_{\cal U} m_\alpha \cdot f,
\]
where
\[
  (g \tensor h) \cdot f := g \ast f \ast h \qquad (f,g,h \in L^1(G))
\] 
and the limit is taken in the norm topology; 
since multiplication by
any element of $L^1(G)$ is compact, this limit does indeed exist. It is
clear that ${\cal Q}$ is a quasi-expectation. However, there are compact groups
$G$ such that $Z(L^1(G))$ is not amenable (\cite{Ste}).
\end{example}
\renewcommand{\baselinestretch}{1.0}
\dated
\vfill
\renewcommand{\baselinestretch}{1.2}
\begin{tabbing} 
{\it Address\/}: \= Department of Mathematical Sciences \\
                 \> University of Alberta \\
                 \> Edmonton \\
                 \> Alberta, T6G 2G1 \\
                 \> Canada \\[\medskipamount]
{\it E-mail\/}:  \> {\tt runde@math.ualberta.ca} \\
                 \> {\tt vrunde@ualberta.ca} \\[\medskipamount]
{\it URL\/}:     \> {\tt http://www.math.ualberta.ca/$^\sim$runde/runde.html}
\end{tabbing}
\end{document}

%% file: format.tex
\renewcommand{\baselinestretch}{1.2}
\newcommand{\dated}{\mbox{} \hfill {\small [{\tt \today}]}}

%% file: mathdefs.tex
\usepackage{amsmath,amssymb,amscd}
\newcommand{\pf}[1]{\trivlist \item[\hskip\labelsep\it #1\ ]}
\newcommand{\varpf}[1]{\trivlist \item[\hskip\labelsep\sc #1:]}
\newcommand{\qedbox}{$\rlap{$\sqcap$}\sqcup$}
\newcommand{\qed}{\qquad \qedbox \endtrivlist}
\newcommand{\varqed}{\hfill \rule{0.6em}{0.6em} \endtrivlist}
\newenvironment{proof}{\pf{Proof}}{\qed}

\newenvironment{remark}{\pf{Remark}}{\endtrivlist}
\newenvironment{remarks}{\pf{Remarks} 
   \begin{enumerate}}{\end{enumerate} \endtrivlist}
\newenvironment{example}{\pf{Example}}{\endtrivlist}
\newenvironment{examples}{\pf{Examples} 
   \begin{enumerate}}{\end{enumerate} \endtrivlist}
\newenvironment{items}{
  \begin{enumerate} 
                    
  }{\end{enumerate}}

\newenvironment{classification}{\noindent\small 2000 {\it Mathematics Subject
Classification\/}:}{\vskip 12pt}

\newcommand{\comps}{{\mathbb C}}
\newcommand{\reals}{{\mathbb R}}

\newcommand{\posints}{{\mathbb N}}

\newcommand{\free}{{\mathbb F}}

\newcommand{\mat}{{\mathbb M}}

\newcommand{\tensor}{\otimes}

\newcommand{\Tensor}{\hat{\otimes}}

\newcommand{\ptensor}{\hat{\otimes}}
\newcommand{\wTensor}{\check{\otimes}}

\newcommand{\cstar}{{C^\ast}}

\newcommand{\id}{{\mathrm{id}}}

\newcommand{\ad}{{\mathrm{ad}}}

\newcommand{\A}{{\mathfrak A}}
\newcommand{\B}{{\mathfrak B}}
\newcommand{\M}{{\mathfrak M}}
\newcommand{\N}{{\mathfrak N}}

%% file: thmdefs.tex
\newtheorem{theorem}{Theorem}[section]
\newtheorem{lemma}[theorem]{Lemma}
\newtheorem{corollary}[theorem]{Corollary}
\newtheorem{proposition}[theorem]{Proposition}
\newtheorem{df}[theorem]{Definition}
\newenvironment{definition}{\begin{df} \rm}{\end{df}}